\documentclass[12pt]{amsart}
\usepackage{graphicx}
\usepackage{subfigure}

\begin{document}

\title[Ocean surface radial velocity imaging]{Ocean surface radial velocity imaging in the AT-INSAR Velocity Bunching Model. A functional approach}

\author{Fabricio Perez, Miguel Angel Moreles, Hector Morales}
\address{\newline Centro de Investigaci\'{o}n en Matem\'{a}ticas\newline
Jalisco s/n, Valenciana\newline
Guanajuato, GTO 36240,Mexico\newline
\textit{email: }fabricio.perez, moreles@cimat.mx
}

\address{\newline Mathematics Deparment\newline
Universidad Autonoma Metropolitana\newline
Ciudad de Mexico,Mexico\newline
\textit{email: }jhmb@xanum.uam.mx
}
\date{} 


\begin{abstract}
This work is concerned with the estimation of radial velocities of sea surface elevations. The data, is a noisy along-track interferometric synthetic aperture radar (AT-INSAR) image. We assume the Velocity Bunching Model. This model relates the complex AT-INSAR image at a point in the image plane, with the radial velocity of a scatterer point in the sea surface. The relation is by means of a nonlinear integral operator mapping radial velocities into AT-INSAR images. Consequently, the estimation of radial velocities,  amounts to the solution of nonlinear integral equations. Our proposal is to solve the latter by Newton's methods on function spaces, the optimize then discretize approach. We show that this continuous version is accurate, and faster than the classical discretize then optimize version. Also a physical comparison is carried out with the interferometric velocities.
\end{abstract}

\maketitle
\tableofcontents

\section{Introduction}

In recent decades, imaging of the surface of our planet Earth has increased with the
appearance and improvement of  tele-detection systems, such as Synthetic Aperture Radar
(SAR).  A SAR system is capable of constructing an image from the information of electromagnetic waves,
which are firstly emitted by the radar and then backscattered by the observed region. 
See Moreira et al (2013) for a review on the subject.

In the case of the ocean surface, Goldstein \& Zebker (1987) developed  the so called, airborne along-track interferometric synthetic aperture radar (AT-INSAR). It has been applied to measure ocean surface currents, ocean surface waves, etc.

Consequently, the understanding of ocean-radar interaction, is of great relevance both in theory and in practice. On the theoretical side, of interest is to derive models of such an interaction. One of such models is the focus of this work, namely, the AT-INSAR Velocity Bunching Model presented in Bao, Br\"{u}ning and Alpers (1997).

This model relates the complex AT-INSAR image at a point in the image plane, with the radial velocity (line of sight velocity component) of a scatterer point in the sea surface.  In the mathematical jargon, this relation is by means of an integral operator mapping radial velocities into AT-INSAR images. If the radial velocity is known in a sea surface region, the AT-INSAR image is readily obtain by quadrature. This is called the direct problem. The purpose of this work, is to  consider the following inverse (imaging) problem:

Given a noisy AT-INSAR image of an unknown scalar field of sea surface elevations, estimate the scalar field of radial velocities of the sea surface elevations.

 A solution to this inverse problem is already contained in Goldstein \& Zebker (1987). The radar data is acquired by two antennas, the fore and aft, carried by a flying platform in the along-track direction at a given velocity. While the aft antenna transmits the signals, both antennas receive the backscattered signals. These are processed separately, then combined interferometrically. It follows that the phase difference caused by the motion of the surface, is proportional to the interferometric velocity. The latter is an approximation to the radial velocity, hence, yields a solution to the inverse problem.
 
The interferometric velocity is used in Hwang et al (2013), for observations of wave breaking in swell-dominant conditions. For further motivation on the imaging problem, see references therein.

A purpose of this work on solving the inverse problem, is to compare the estimated radial velocity field with that of the interferometric velocities. It serves as a query on the At-INSAR Velocity Bunching model. 

As a first study,  we consider an AT-INSAR image generated by a a swell sea. Then, we estimate radial velocities by solving the AT-INSAR Velocity Bunching integral equation.  

The core of the paper is on developing numerical methods for solving the underlying integral equation. It will 
become apparent that the latter is nonlinear and oscillatory, which makes the solution challenging. 

A classical approach is to discretize the integral equation and solve the resulting nonlinear system by Newton's methods. The so called discretize then optimize approach. 

Alternatively we opt to postpone discretization \emph{until the last minute}, that is, the optimize then discretize approach. It is proven to be more efficient, Stuart (2010), and sometimes necessary, Zuazua (2005).

We develop two modified Newton's methods on function spaces for solution. First a nonlinear system, second
as a nonlinear least squares problem. Derivatives are computed in the sense of Fr\'{e}chet. See for instance Cheney (2001), for the required Functional Analysis. 

Both solutions are mutually cross-validating. For comparison we implement also a discretize then optimize approach,
in the case of minimization. As expected, the former perform better and on execution time are considerably faster.

A physical comparison is also carried out between the bulk kinetic energy on the ocean surface area under study, associated to the estimated fiel of radial velocities and that of the interferometric velocities. The comparison is in terms of relative errors, again, the latter is outperformed.

We work with synthetic data, we generate a field of surface elevations following the classical variance spectra to surfaces approach.  We obtain a random 2-D realization of a sea surface. Following Mobley (2016), we develop our own implementation.

\section{Materials and Methods}

In this section we pose the imaging problem of interest and the modified Newton's methods for solution. We follow a functional approach, that is, we optimize on normed vector spaces of functions. 

More precisely, we shall consider all function spaces as subspaces of $L^2\big((a,b)\big)$, the space of square summable complex functions.  
For two such functions, $\phi$, $\psi$ the inner product is
\begin{equation} \label{eq:L2_innerProduct}
\big\langle \phi \, , \psi \big\rangle \quad = \quad \int_a^b \! \phi(x)\,\overline{\psi(x)}\,dx
\end{equation}
We shall use freely all well known hilbertian properties of $L^2$, see for instance Cheney (2001).

\subsection{Problem statement}

The point of departure is the AT-INSAR configuration as in Goldstein \& Zebker (1987).

 The AT-INSAR image is acquired by two antennas, the fore and aft, carried by a flying platform in the along-track direction at velocity $V$. The antennas are separated by a a $2B$ distance. We assume that the system operates in mode 1, the aft antenna transmits radar signals,
and both antennas receive the backscattered signals. The wavenumber of the incoming electromagnetic wave is denoted by $k_r$.

Let $\mathbf{x}=(x,y)$ be in the reference frame for the sea surface $z(\mathbf{x})$. $x$ is the coordinate in ground range (cross-track), and $y$ the coordinate in azimuth (along-track).

Let $I(\mathbf{x}_R)$ be the AT-NSAR image at the position $\mathbf{x}_R=(x_R,y_R)$ that is associated with the scatterer $P_\mathbf{x}$ at the point $(\mathbf{x}, z(\mathbf{x}))$. The distance from the median of the two antennas to the point $(\mathbf{x}, 0)$, is denoted by $R$. Also, 
denote by $\tau_s$, the scene coherence time.

We assume the AT-INSAR Velocity Bunching Model for a complex AT-INSAR image $I_{vb}$ as introduced in Bao, Bruning and Alpers (1997).

In this model,  the AT-INSAR (single-look) integration time, $T_0$ is regarded small compared to the period of the dominant ocean wave.  Hence, the normalized radar cross section (NRCS)  $\sigma(\mathbf{x},t)$ and the radial velocity $u_r(\mathbf{x},t)$ vary little and are approximated by quantities independent of time, denoted by  $\sigma_0(\mathbf{x})$, $u_r(\mathbf{x})$ respectively. In particular, for the radial velocity a first order approximation in time iabout $t_0=\mathbf{x}/V$ is used,
\[
u_r(\mathbf{x},t) \approx u_r(\mathbf{x})+a_r(\mathbf{x})(t-t_0).
\]
Here $a_r(\mathbf{x})$ is the radial acceleration.

 With these considerations, the expression for $I_{vb}(y_R)\equiv I_{vb}(\mathbf{x}_R)$ is,

\begin{eqnarray}  
I_{vb}(y_R) &=& \frac{\pi T_0^2 \rho_a}{2} exp\bigg[-\frac{4B^2}{V^2T_0^2}\bigg] \ \int_{-\infty}^{+\infty} \frac{\sigma_0(\mathbf{x})}{\rho_a'(\mathbf{x})} \nonumber \\
&& \times \ exp\bigg[-2jk_r \frac{B}{V} u_r(\mathbf{x})\bigg]  exp\bigg[\frac{4B^2\rho_a^2}{V^2\,T_0^2\,\rho_a'^2(\mathbf{x})}\bigg] \nonumber \\
&& \times \ exp \bigg[ \frac{2jBk_r}{R}\bigg(\frac{2\rho_a^2}{\rho_a'^2(\mathbf{x})} - 1\bigg) \bigg(y_R - y - \frac{R}{V}u_r(\mathbf{x})\bigg) \bigg] \nonumber \\
&& \times \ exp \bigg[ -\frac{\pi^2}{\rho_a'^2(\mathbf{x})} \bigg(y_R - y - \frac{R}{V}u_r(\mathbf{x})\bigg)^2 \bigg] dy \nonumber \\ \nonumber \\
\rho_a'(\mathbf{x})  & = & \Bigg\{\rho_a^2 + \bigg[\frac{\pi}{2}\,\frac{T_0R}{V}\,a_r(\mathbf{x})\bigg]^2 + \frac{\rho_a^2T_0^2}{\tau_s^2}\Bigg\}^{1/2} \nonumber
\end{eqnarray}

Here $\rho_a'(\mathbf{x})$ denotes the degraded single-look azimuthal resolution, and $\rho_a=\lambda_rR/(2VT_0)$ is the full-bandwidth,
single-look azimuthal resolution for stationary targets, where $\lambda_r$ denotes the radar wavelength.

The inverse (imaging) problem of interest is: Given AT-INSAR noisy data $\text{D}$ of an unknown sea surface $z$, and given the radar parameters, estimate the radial velocities $u_\text{r}$ of $z$. 

It is assumed that $D(\mathbf{x}_R)$ is a complex AT-INSAR image $I_{vb}(\mathbf{x}_R)$ corrupted by additive noise $\eta$, namely $D=I_{vb}+\eta$. 

Notice that the problem amounts to solving an oscillatory nonlinear integral equation for $u_r(\cdot)\equiv u_r(x,\cdot)$, for each fixed $x$ in the cross-track coordinate in the observation area.

For later reference, let us define
\begin{equation}
A = \frac{\pi T_0^2 \rho_a}{2} \text{exp} \bigg[-\frac{4B^2}{V^2T_0^2}\bigg] 
\end{equation}

Denoting the integrand by $f_\text{vb}$, a scalar complex-valued function, we have.
\begin{equation}
I_{vb}(y_R) =  A\int_{-\infty}^{+\infty}\!\!\!f_\text{vb}(u_r(\mathbf{x}),a_r(\mathbf{x}),\mathbf{x},y_R)\, dy
\end{equation}

\bigskip

\subsubsection*{Remark}
In our case study, it is found a fortiori, that variations of this integral operator with respect to the radial acceleration are negligible.
Consequently, the terms involving the latter in the Newton's methods that follow, are discarded. It is possible to show this mathematically,
but we focus on the numerical results.

\subsection{Newton's method for the nonlinear integral equation}

The nonlinear integral equation above, implicitly defines a map between some function spaces $\mathcal{V}$ and $\mathcal{W}$. Namely,
\[
\mathcal{V}\to\mathcal{W},\quad u_r\mapsto I_{vb}(u_r).
\]

To pose the inversion problem as the solution of a nonlinear integral equation, consider the residual map 
\[
\mathcal{F}: \mathcal{V} \rightarrow \mathcal{W},  \quad \mathcal{F}(u_r)=D-I_{vb}(u_r).
\]

The problem is to find $u_r$ such that
\[
\mathcal{F}(u_r)=0.
\]

Assuming Frechet differentiability, we apply the Newton's method.
 
 Given the initial guess $u_r^0 \in \mathcal{V}$, solve at each iteration $k$ for the function $h$

\begin{equation}
\mathcal{F}'(u_r^k)\,h \ = \ -\mathcal{F}(u_r^k) 
\label{Nwt1}
\end{equation}
and update
\begin{equation}
u_r^{k+1}  = u_r^k \, + \, h.
\label{Nwt2}
\end{equation}
Here $\mathcal{F}'$ is the Frechet derivative of $\mathcal{F}$. It follows that $\mathcal{F}'=-I_{vb}'$, and
 
\[
 I_{vb}'(u_r) h 
= A \int_{-\infty}^\infty \bigg[\frac{\partial f_vb}{\partial u_r}\bigg] h(y)\,dy 
\]
where $(\partial f_{vb}/\partial u_r)$ is the vector calculus derivative of $f_{vb}$ with respect to $u_r$,
\begin{equation}
\frac{\partial f_{vb}}{\partial u_r} =  \bigg[\frac{2\,\pi^2 R \,C}{V \rho_a'\,\!\!^2}  \ - \ j\,\frac{4\,B\,k_r\,\rho_a^2}{V\rho_a'\,\!\!^2}\bigg] \, f_{vb} 
\end{equation}

We remark that the scheme (\ref{Nwt1}),(\ref{Nwt2}) is the continuous (infinite dimensional) version of the Newton's method. 

We discretize with an appropriate quadrature to obtain a finite dimensional residual $F: \mathbb{R}^{N_y} \rightarrow \mathbb{R}^{2 N_y}$.

The Newton's method reads: Given an initial guess $\vec{u}_r\,\!\!^0 \in \mathbb{R}^{N_y}$, solve at each iteration $k$:

\begin{eqnarray}
F'(\vec{u}_r\,\!\!^k)\,\vec{h}^k \ &=& \ -F(\vec{u}_r\,\!\!^k) \label{eq:NewtonsMethod_1_numerical} \\
\vec{u}_r\,\!\!^{k+1} \ &=& \ \vec{u}_r\,\!\!^k \, + \, \vec{h}^k \label{eq:NewtonsMethod_2_numerical}
\end{eqnarray}

The linear system to solve in each iteration is overdetermined and ill conditioned. Thus the numerical solution is somewhat challenging. 

We apply Tikhonov regularization using the SVD decomposition of $F'(\vec{u}_\text{r}\,\!\!^k)$. Namely

\begin{eqnarray} \label{eq:TikhonovSolution_SVD}
\vec{h}^k  \ &\approx& \ \vec{h}_{\alpha_T}^k \quad = \quad - \sum_{i=1}^{N_y}\frac{\sigma_i}{\sigma_i^2 + \alpha_T} \Big[\vec{u}_i\,\!\!\cdot F(\vec{u}_r\,\!\!^k)\Big] \vec{v}_i
\end{eqnarray}
where $\sigma_i \in \mathbb{R}_{\geq 0}$ is the $i$-th largest singular value of $F'(\vec{u}_r\,\!\!^k)$, $\vec{u}_i \in \mathbb{R}^{2 N_y}$ is the corresponding $i$-th left singular vector, and $\vec{v}_i \in \mathbb{R}^{N_y}$ is the corresponding $i$-th right singular vector. 
 
The regularization parameter is set to the square of the first (and largest) singular value of $F'(\vec{u}_r\,\!\!^k)$:

\begin{eqnarray} \label{eq:alpha_Tikhonov}
\alpha_T &\leftarrow& \sigma_1^2
\end{eqnarray}

\subsection{A minimization approach}

For comparison, we develop a descent method for the corresponding nonlinear squares problem.  The underlying 
functional,

\[
 \mathcal{G}: \mathcal{V} \rightarrow [0,\infty),
\]
is given by

\[
\mathcal{G}(u_\text{r})=\frac{1}{2}\Vert 
\mathcal{F}(u_\text{r}) \Vert^2 = \frac{1}{2} \langle \mathcal{F}(u_\text{r}) \, , \, \mathcal{F}(u_\text{r}) \rangle
\]

It is readily seen that the Fr\'echet derivative of $\mathcal{G}$ at $u_\text{r}$ applied to $h$ is:

\[
\mathcal{G}'(u_\text{r})h = \frac{1}{2}\left[ \langle h , ( \mathcal{F}'(u_\text{r}))^*\mathcal{F}(u_\text{r})\rangle + \overline{\langle h , ( \mathcal{F}'(u_\text{r}))^*\mathcal{F}(u_\text{r})\rangle}\right]
\]
where $( \mathcal{F}'(u_\text{r}))^*$ is the adjoint operator of  $\mathcal{F}'(u_\text{r})$.

Hence
\[
\mathcal{G}'(u_\text{r})h =  \langle h , Re \left\lbrace
( \mathcal{F}'(u_\text{r}))^*\mathcal{F}(u_\text{r})
\right\rbrace
\rangle
\]

Substituting $\mathcal{F}'(u_\text{r})$, it follows that

\begin{eqnarray} \label{eq:Gprime_widetilde_ur_h_2nd}
\mathcal{G}'(u_r)h \ &=& \ \Bigg\langle h, Re\,\bigg\{\!- A \int_{-\infty}^\infty \bigg[\frac{\partial f_{vb}}{\partial u_r}\bigg]\,\overline{\bigg[\big[\mathcal{F}(u_r)\big](y_R})\bigg]\,dy_R\, \bigg\} \Bigg\rangle
\end{eqnarray}

By the Riesz representation theorem, the second argument of the inner product is the continuous gradient of $\mathcal{G}$ at $u_r$. 
That is, the function $\nabla \mathcal{G}(u_r)$,  defined as:
\begin{eqnarray}
\nabla \mathcal{G}(u_r) \ &=& Re\,\bigg\{\!- A \int_{-\infty}^\infty \bigg[\frac{\partial f_{vb}}{\partial u_r}\bigg]\,\overline{\bigg[\big[\mathcal{F}(u_r)\big](y_R})\bigg]\,dy_R\, \bigg\}  \label{eq:gradGwidetildeur}
\end{eqnarray}

The minimization problem is solved by the BFGS method with known gradient. In our case,  a discrete version of (\ref{eq:gradGwidetildeur}). 

\section{Synthetic data}

The area under study is a square $Q$, with $1280m$ side length centered at the origin. A $128\times 128$ uniform square mesh is considered, with computation points at the vertices.

The purpose of this section is to specify the parameters to simulate the ocean surface $z: Q \rightarrow \mathbb{R}$ and its associated AT-INSAR data $D: Q_R\rightarrow \mathbb{C}$. 

\subsection{Sea surface}
We follow the classical variance spectra to surfaces approach, to generate a random 2-D realization of a sea surface. See Mobley (2016).

From knowledge from real ocean surfaces one starts  with an omnidirectional spectrum. For a swell sea  we use at frequency $k$ the spectrum,Bao, Bruening and Alpers (1997).

\begin{equation}
\mathcal{S}_S(k) = \frac{\alpha_S}{2\,k^3}\, exp\Bigg[-\frac{5}{4}\bigg(\frac{k}{k_S}\bigg)^{-2}\Bigg] \,\,\gamma_S\,\!^{G_S} \nonumber
\end{equation}
where
 \begin{equation} 
 G_S = exp\Bigg[\!\!-\frac{1}{2}\,\frac{\big(k^{1/2} \ - \ k_S^{1/2}\big)^2}{\sigma_S^2\,k_S} \Bigg] \nonumber
\end{equation}

Here,
\begin{itemize}

\item $\alpha_S=0.212 \times 10^{-3}$, the energy scale of $\mathcal{S}_S$.  

\item $k_S \ = \ 2\pi \,(\lambda_S)^{-1}$,  the spatial peak frequency of $S_S$. 

\item $\gamma_S = 10$, the peak enhancement factor of \! $\mathcal{S}_S$.

\item $\sigma_\text{S}$, the spectral width centered at $k_\text{S}$,
\end{itemize}

\begin{equation}
\sigma_\text{S} \quad =  \quad
\begin{cases}
0.07 & \quad \text{for} \ \  k \, \leq \, k_\text{S} \\
0.09 & \quad \text{for} \ \  k \, > \, k_\text{S}
\end{cases}
\end{equation}

Next a spreading function is used. In this case, the two-sided cosine-power model,
\[
\Phi_\text{cp2}(k,\phi) \ = \ \frac{1}{2} \, N_p \, \vert\text{cos}\,(\phi - \phi_w) \vert^{2p}, \qquad \vert\phi - \phi_w\vert\, \leq \, \pi,
\]
Leading to the two-sided directional swell spectrum,
\[
\mathcal{\hat{S}}_{2C}(k_x, k_y) \ = \ \frac{1}{k} \,\mathcal{S}_1(k) \,\Phi_\text{cp2}(k,\phi).
\]

Where $(k,\phi)$ \! and \! $(k_x, k_y)$ \! are the equivalent polar and cartesian coordinates, respectively.

We obtain a particular instance of the ocean variance spectrum $\hat{z}$. The ocean surface $\text{z}$ is obtained by computing the discrete
inverse Fourier Transform of $\hat{z}$. 

We have developed our own software for ocean surface simulation. Let us show some graphics from omnidirectional swell spectrum to ocean and radial velocity surfaces.

\begin{figure}[!htb]
\centering
\subfigure[Omnidirectional spectrum $\mathcal{S}_\text{S}$: continuous (line) and sampled (dots).]{\label{fig:A07_S_1S}
\includegraphics[width=5.8in]{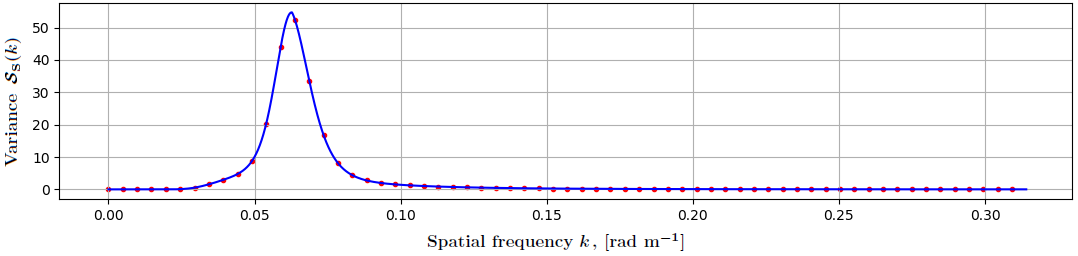}} 
\vspace{0.2cm}
\subfigure[Sampled directional spectrum $\hat{\text{S}}_\text{2C}$.]{\label{fig:A02_Psi_2S}
\includegraphics[width=2.67in]{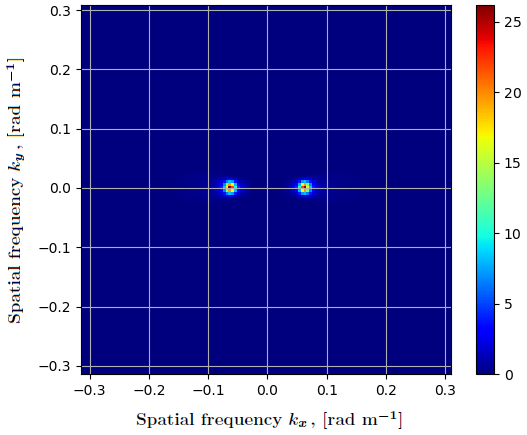}} 
\subfigure[Sampled-variance spectrum $\vert\hat{\text{z}}\vert^2$.]{\label{fig:A05_Var}
\includegraphics[width=2.82in]{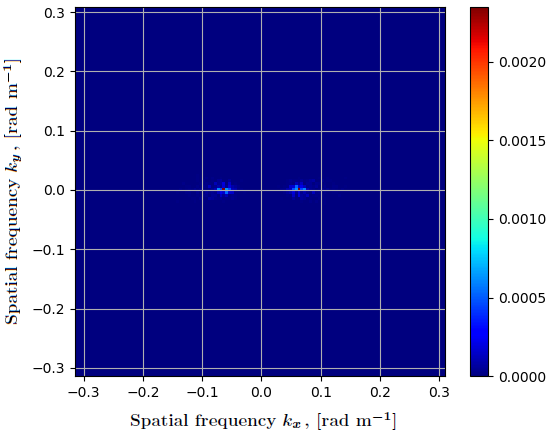}} 
\caption[]{Ocean spectra to construct and validate the simulated ocean surface.}
\label{fig:D_and_derived_fields}
\end{figure}

In the frequency domain, figures \ref{fig:A07_S_1S}, \ref{fig:A02_Psi_2S}, \ref{fig:A05_Var} depict the omnidirectional spectrum 
$\mathcal{S}_\text{S}$, the sampled directional spectrum $\hat{\text{S}}_\text{2C}$ and the sampled-variance spectrum $\vert\hat{\text{z}}\vert^2$, respectively.

\begin{figure}[!htbp]
	\centering
		\includegraphics[scale=0.59]{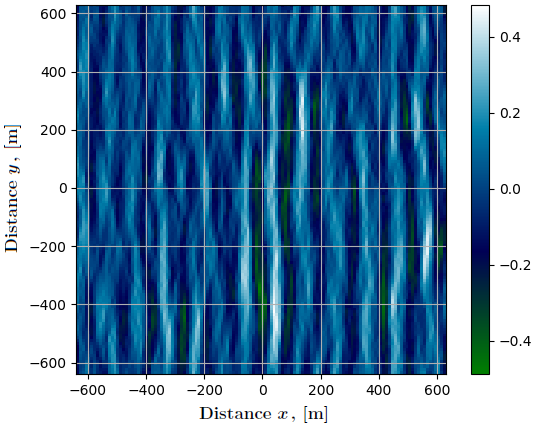} 
	\caption{The simulated ocean surface $\text{z}$.}
	\label{fig:A06_SurfaceElevation_z}
\end{figure}

The ocean surface $\text{z}$ is shown in figure \ref{fig:A06_SurfaceElevation_z}, where three properties of such surface can be visualised: there is a regular pattern of waves whose directions are very close to the wind direction $\phi_w = 0$ [rad], the majority of wavelengths are around $\lambda_\text{S} = 100$ [m], and a big amount of measured wave heights are well characterised by $H_{m_0} = 0.586215$ [m].

\begin{figure}[!htbp]
	\centering
		\includegraphics[width=2.92in]{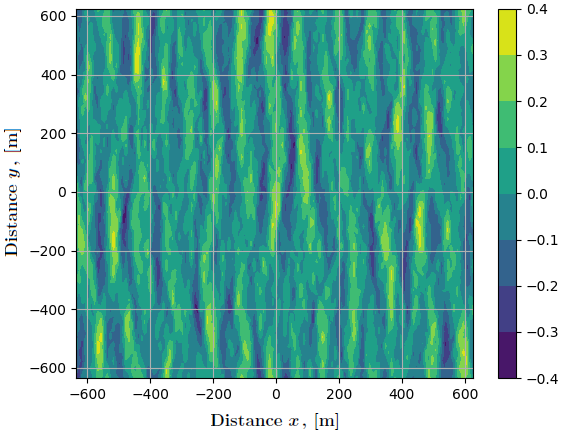}
	\caption{Scalar fields of radial velocities.}
	\label{fig:B02_0_ur}
\end{figure}

Finally, the scalar field of radial velocities in Figure \ref{fig:B02_0_ur}

\subsection{The AT-INSAR data} 

The configuration of the AT-INSAR system is that of Bao, Bruening and Alpers (1997). For each pointt $\mathbf{x}$  in the grid, $I_{vb}(\mathbf{x})$ is approximated by quadrature. Then noise $\eta(\mathbf{x})$ is added. More precisely,

 \[
 \eta(\mathbf{x}) = 2^{-1/2}\big[a_\eta(\mathbf{x}) \, + \, jb_\eta(\mathbf{x})\big],
 \] 
 where elements $a_\eta(\mathbf{x})$ and $b_\eta(\mathbf{x})$ are independent real Gaussian random variables with mean $0$ and variance 
 $\sigma_\eta^2$. 
 
The standard deviation $\sigma_\eta$, is taken as:
\begin{equation}
\sigma_\eta = \left[
10^{SNR/20}
\right]^{-1},
\end{equation}
where $SNR$ is the signal to noise ratio in [dB]. Below, we report results for the value $SNR = 174$ [dB]. 

\subsection{Interferometric velocities}

Given the AT-INSAR noisy image $D$, the interefometric phase $\Phi_{ATI}$ is given by $\Phi_{ATI} = \angle D$. It is proportional to the inetrferometric velocity $u_{ATI}$, Goldstein \& Zebker (1987). Namely
\[
 u_{ATI}(\mathbf{x})=-\frac{\lambda_r}{4\pi}\,\frac{V}{B}\,\Phi_{ATI}(\mathbf{x}_R).
\]

In applications, $u_{ATI}$ is used as an approximation of the radial velocity $u_r$. We gauge this approximation in the results that follow.

\section{Numerical Results}

\subsection{Radial Velocity Imaging}

For each $x$ in the cross-track coordinate we solve both, the Nonlinear System (NL), and the Functional Minimization (FM) associated to the 
nonlinear integral equation. In total there are $128$ problems ordered from left to right, from $0$ to $127$. 

We stress that we follow the optimize then discretize approach. For comparison, we use the discretize then optimize in the minimization problem.
Results with the latter shall be referred as DFM. 

The main difference is the computation of the gradient, as usual,  it is approximated with appropriate finite differences of the discretized (finite dimensional) functional.

To avoid bias, we start all iterative methods with $u_r^0\equiv 0$ as initial guess. Results are remarkable homogeneous for the 128 problems. 

First we show the fit for problems $0$ and $64$ in figures 4 and 5, respectively. 
The interferometric velocity is included.

\begin{figure}[htbp]
	\centering
		\includegraphics[scale=0.23]{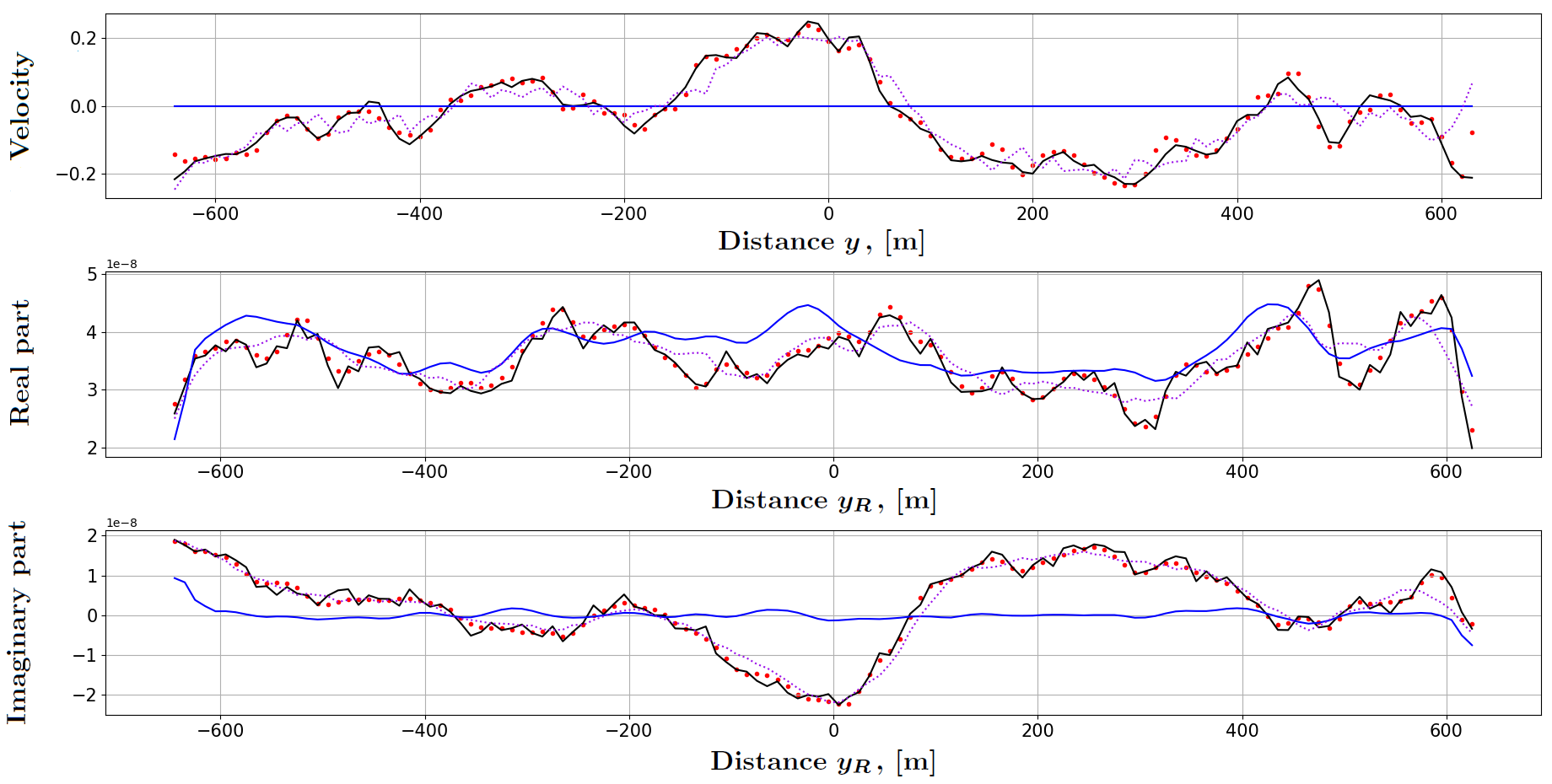}     
           \caption[]{\scriptsize TOP: radial velocity $\text{u}_\text{r}$ (black), initial point $\text{u}_\text{r}\,\!\!^0$ (blue), interferometric velocity $\text{u}_\text{ATI}$ (violet), estimated radial velocity $\text{u}_\text{r}\,\!\!^*$ (red); MIDDLE: $\text{Re}\big\{\text{D}\big\}$ (black), $\text{Re}\big\{\,\text{I}_0\big\}$ (blue), $\text{Re}\big\{\,\text{I}_\text{ATI}\big\}$ (violet), $\text{Re}\big\{\,\text{I}_*\big\}$ (red); BOTTOM: $\text{Im}\big\{\text{D}\big\}$ (black), $\text{Im}\big\{\,\text{I}_0\big\}$ (blue), $\text{Im}\big\{\,\text{I}_\text{ATI}\big\}$ (violet), $\text{Im}\big\{\,\text{I}_*\big\}$ (red).}
           \label{fig:ssnle_p0}  
\end{figure}

\begin{figure}[htbp]
	\centering
		\includegraphics[scale=0.23]{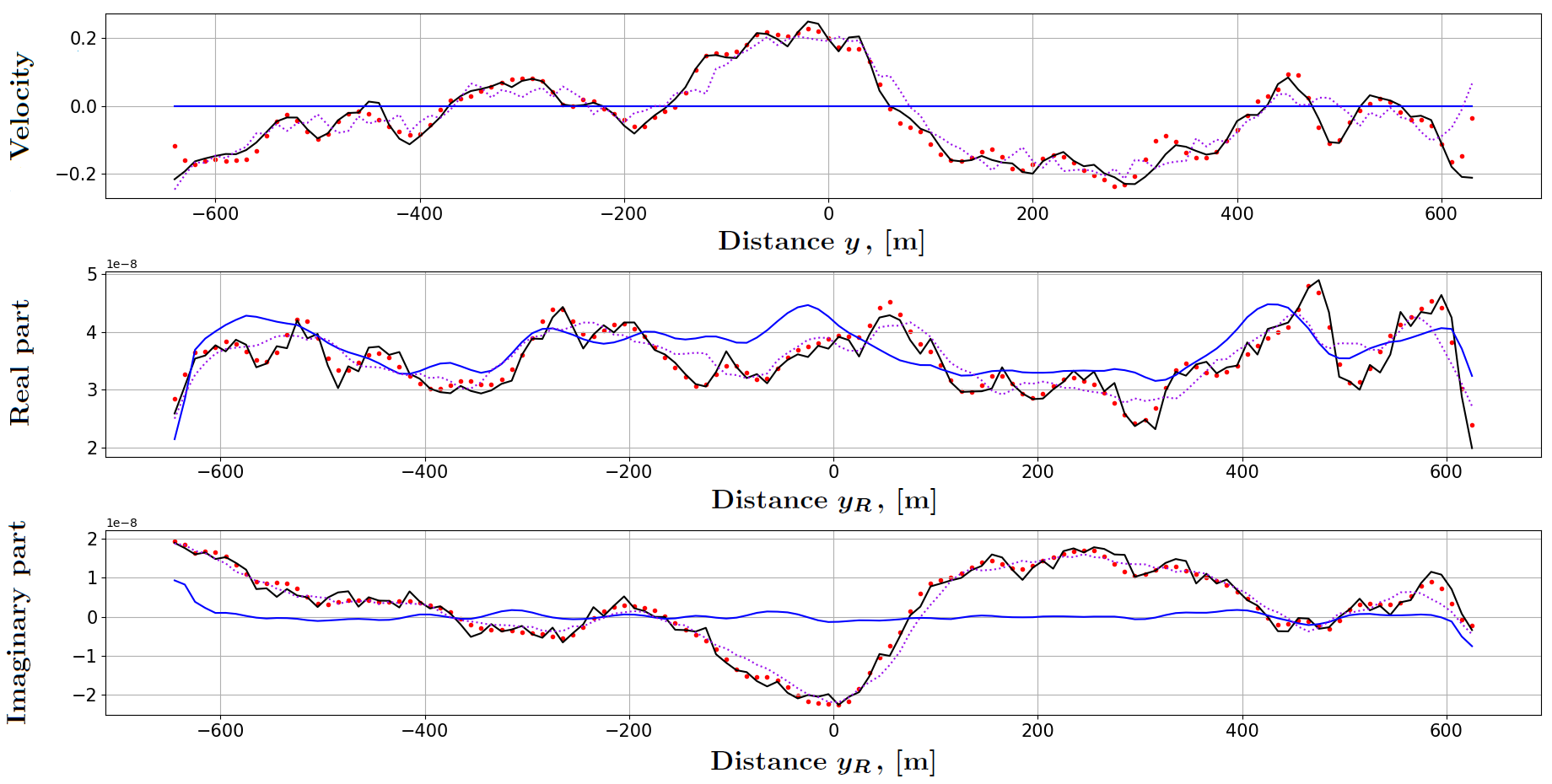}
           \caption[]{TOP: radial velocity $\text{u}_\text{r}$ (black), initial point $\text{u}_\text{r}\,\!\!^0$ (blue), interferometric velocity $\text{u}_\text{ATI}$ (violet), estimated radial velocity $\text{u}_\text{r}\,\!\!^*$ (red); MIDDLE: $\text{Re}\big\{\text{D}\big\}$ (black), $\text{Re}\big\{\,\text{I}_0\big\}$ (blue), $\text{Re}\big\{\,\text{I}_\text{ATI}\big\}$ (violet), $\text{Re}\big\{\,\text{I}_*\big\}$ (red); BOTTOM: $\text{Im}\big\{\text{D}\big\}$ (black), $\text{Im}\big\{\,\text{I}_0\big\}$ (blue), $\text{Im}\big\{\,\text{I}_\text{ATI}\big\}$ (violet), $\text{Im}\big\{\,\text{I}_*\big\}$ (red).}
           \label{fig:umf_p0}
\end{figure}

The corresponding fitting results for problem 64 are shown in figures 6 and 7

\begin{figure}[htbp]
	\centering
		\includegraphics[scale=0.23]{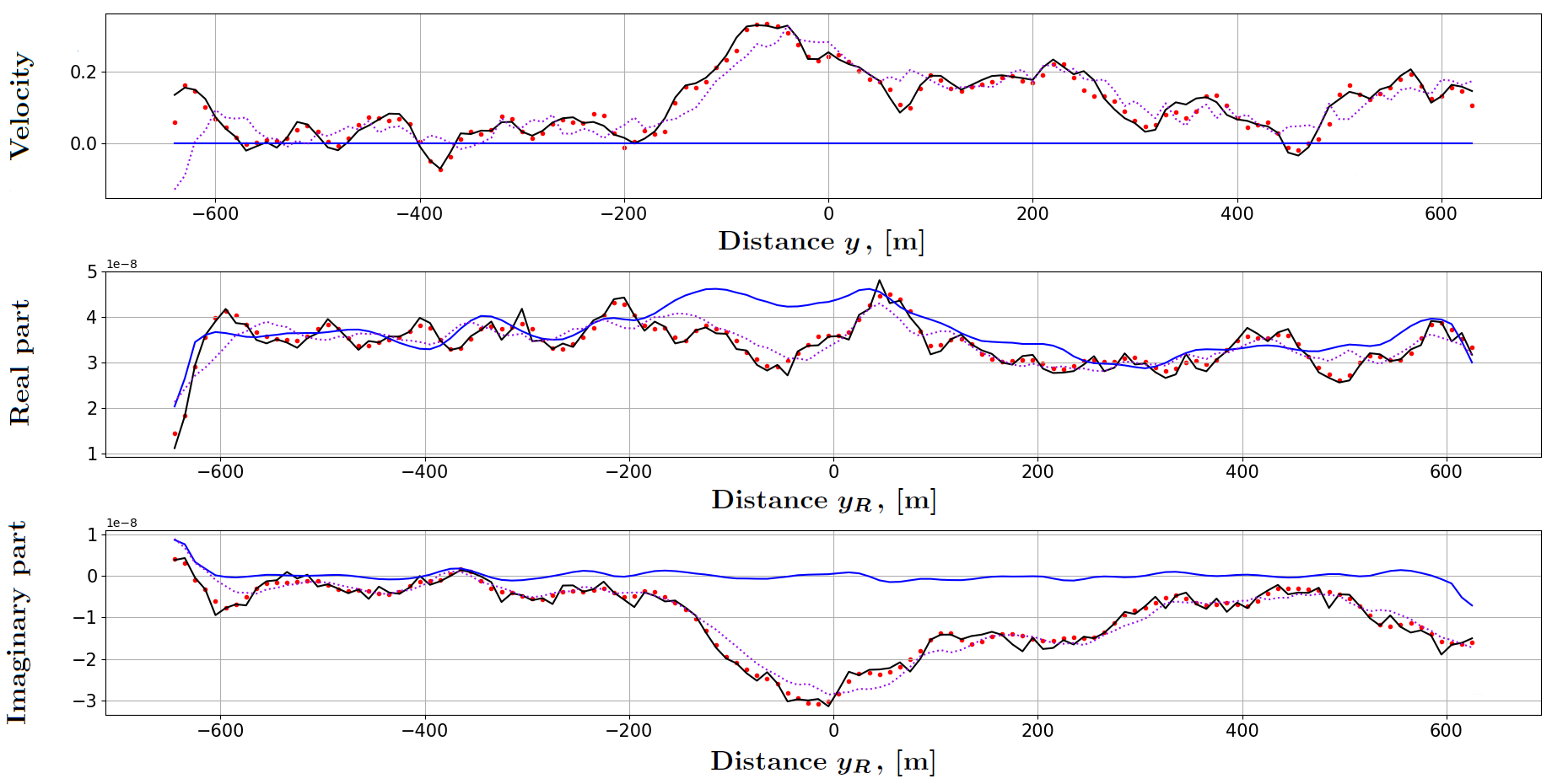}
	\label{fig:ssnle_p64}
           \caption[]{ TOP: radial velocity $\text{u}_\text{r}$ (black), initial point $\text{u}_\text{r}\,\!\!^0$ (blue), interferometric velocity $\text{u}_\text{ATI}$ (violet), estimated radial velocity $\text{u}_\text{r}\,\!\!^*$ (red); MIDDLE: $\text{Re}\big\{\text{D}\big\}$ (black), $\text{Re}\big\{\,\text{I}_0\big\}$ (blue), $\text{Re}\big\{\,\text{I}_\text{ATI}\big\}$ (violet), $\text{Re}\big\{\,\text{I}_*\big\}$ (red); BOTTOM: $\text{Im}\big\{\text{D}\big\}$ (black), $\text{Im}\big\{\,\text{I}_0\big\}$ (blue), $\text{Im}\big\{\,\text{I}_\text{ATI}\big\}$ (violet), $\text{Im}\big\{\,\text{I}_*\big\}$ (red).}
\end{figure}

\begin{figure}[htbp]
	\centering
		\includegraphics[scale=0.23]{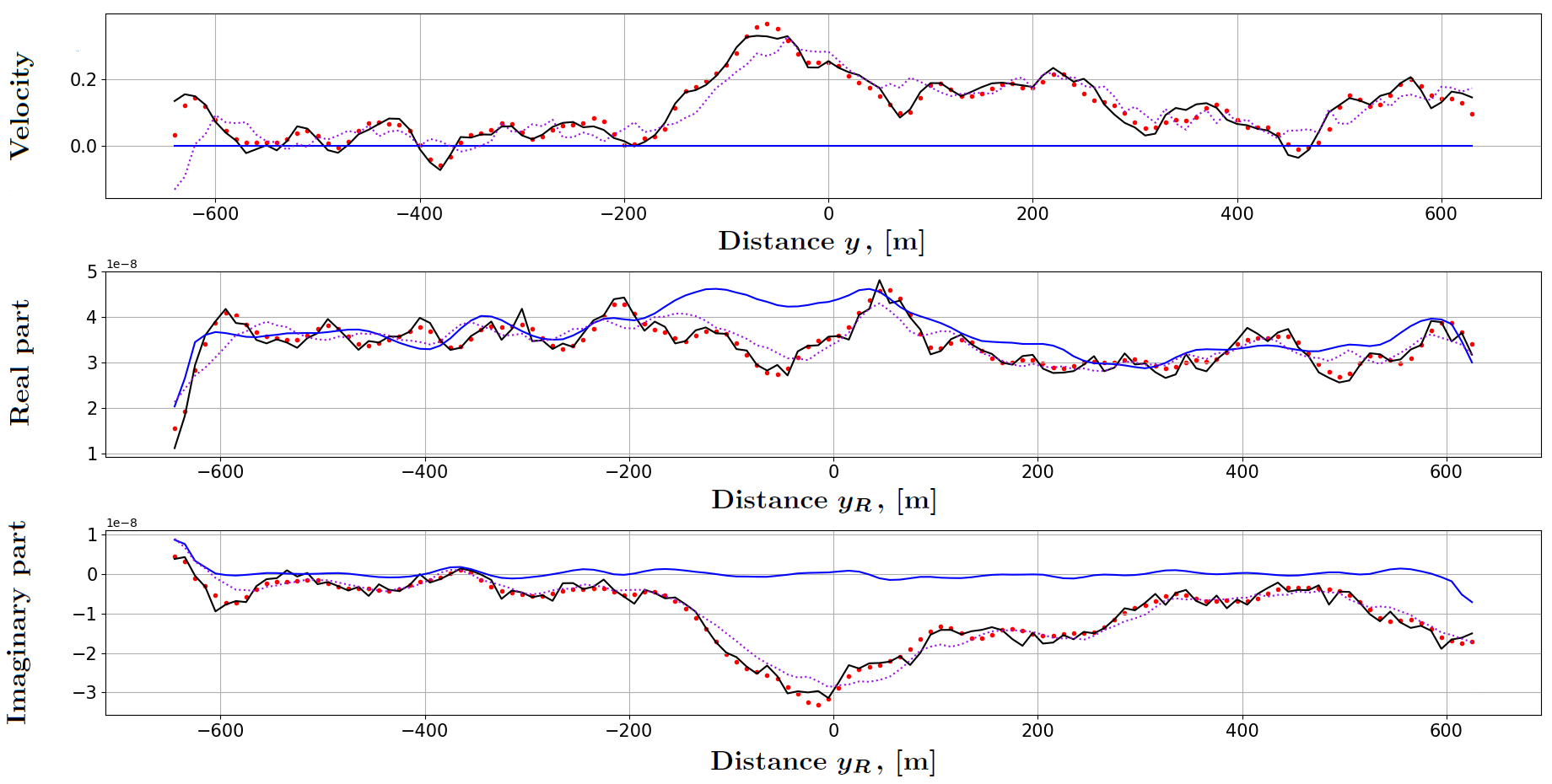}
	\label{fig:umf_p64}
           \caption[]{ TOP: radial velocity $\text{u}_\text{r}$ (black), initial point $\text{u}_\text{r}\,\!\!^0$ (blue), interferometric velocity $\text{u}_\text{ATI}$ (violet), estimated radial velocity $\text{u}_\text{r}\,\!\!^*$ (red); MIDDLE: $\text{Re}\big\{\text{D}\big\}$ (black), $\text{Re}\big\{\,\text{I}_0\big\}$ (blue), $\text{Re}\big\{\,\text{I}_\text{ATI}\big\}$ (violet), $\text{Re}\big\{\,\text{I}_*\big\}$ (red); BOTTOM: $\text{Im}\big\{\text{D}\big\}$ (black), $\text{Im}\big\{\,\text{I}_0\big\}$ (blue), $\text{Im}\big\{\,\text{I}_\text{ATI}\big\}$ (violet), $\text{Im}\big\{\,\text{I}_*\big\}$ (red).}
\end{figure}

In figure 8 and 9, we show the RMSE of the nonlinear system solution $u_r^*$--NL,  the Discrete Functional
Minimization solution $u_r^*$--DFM, and the interferometric velocity solution $u_{ATI}$. Each point in the horizontal axis corresponds
to a fitting problem, 128 in total.

It is apparent that in both cases the functional (infinite dimensional) approach performs better that the DFM (finite dimensional) solution.

\begin{figure}[!htbp]
	\centering
		\includegraphics[scale=0.33]{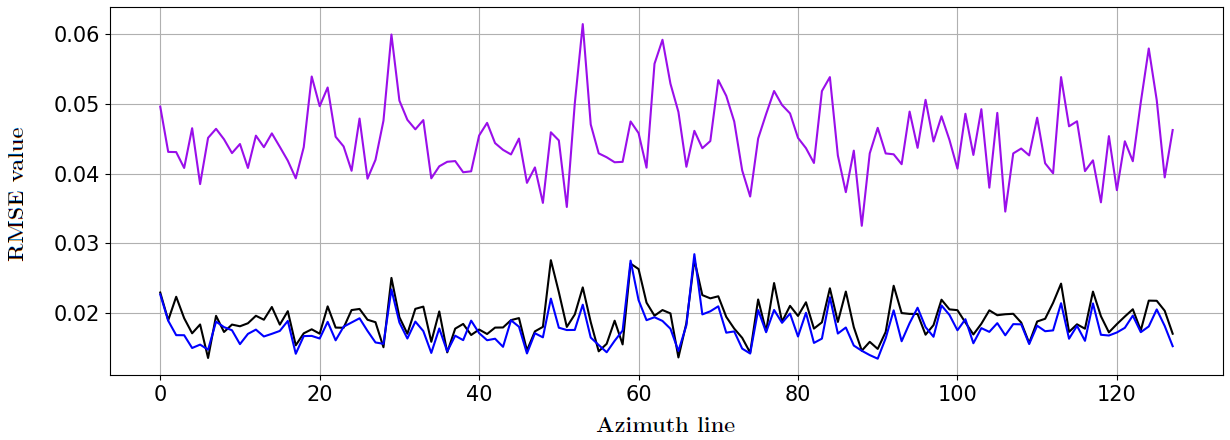}
		\label{fig:NL_vs_DFM}
		\caption{Blue: $\widetilde{\text{u}}_\text{r}^*$--NL \quad Red: $\widetilde{\text{u}}_\text{r}^*$--FM \quad Black: $\widetilde{\text{u}}_\text{r}^*$--DFM \quad Violet: $\text{u}_\text{ATI}$}
\end{figure}

\begin{figure}[!htbp]
	\centering
		\includegraphics[scale=0.33]{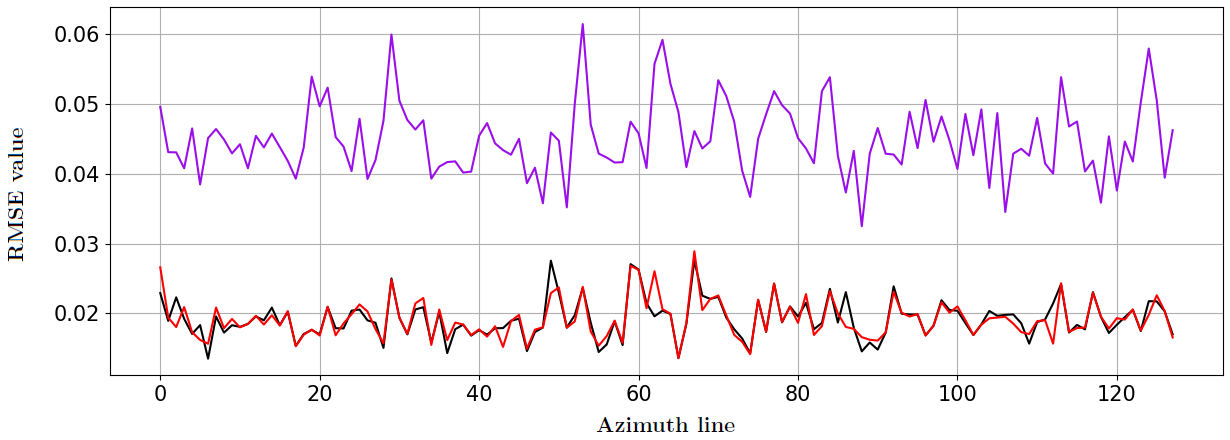}
		\label{fig:FM_vs_DFM}
		\caption{Blue: $\widetilde{\text{u}}_\text{r}^*$--NL \quad Red: $\widetilde{\text{u}}_\text{r}^*$--FM \quad Black: $\widetilde{\text{u}}_\text{r}^*$--DFM \quad Violet: $\text{u}_\text{ATI}$}
\end{figure}

\subsection{Computational Efficiency}

In practice, imaging problems are computationally expensive. In the discretize then optimize approach, the approximation of derivatives by finite 
differences is costly. Having the exact derivative, and postponing discretization until \emph{the last minute} is in general more efficient.

In figures 10 and 11, we show execution times for the 128 inversion problems. 

It is noticed that the functional versions are at least three orders of magnitud faster.

\begin{figure}[!htbp]
	\centering
		\includegraphics[scale=0.33]{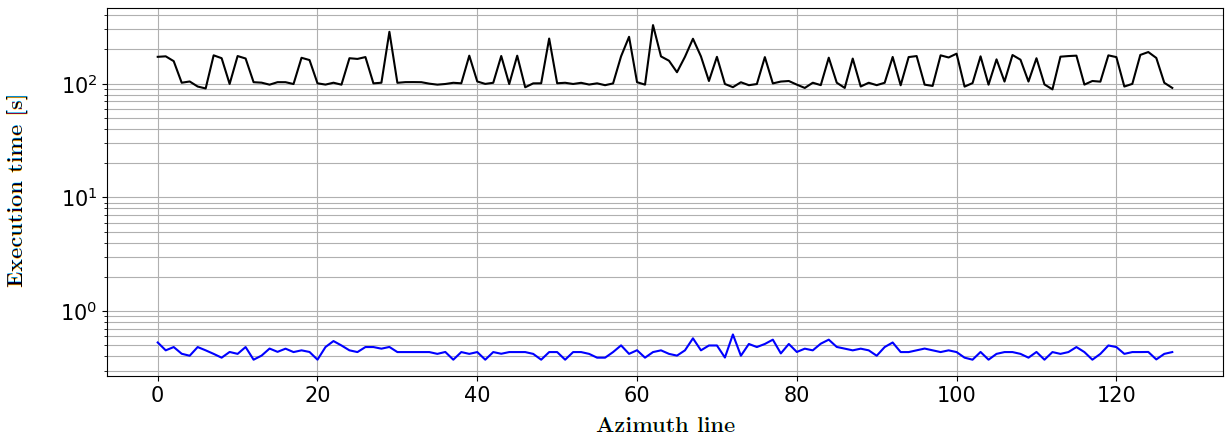}
		\label{fig:NL_vs_DFM_exT}
		\caption{Blue: $\widetilde{\text{u}}_\text{r}^*$--NL \quad Red: $\widetilde{\text{u}}_\text{r}^*$--FM \quad Black: $\widetilde{\text{u}}_\text{r}^*$--DFM \quad Violet: $\text{u}_\text{ATI}$}
\end{figure}

\begin{figure}[!htbp]
	\centering
		\includegraphics[scale=0.33]{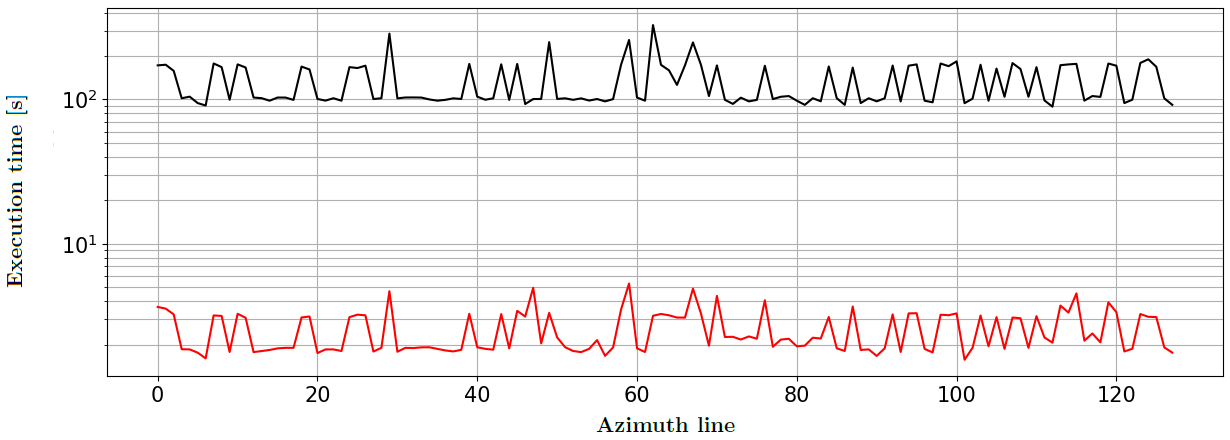}
		\label{fig:FM_vs_DFM_exT}
		\caption{Blue: $\widetilde{\text{u}}_\text{r}^*$--NL \qquad Red: $\widetilde{\text{u}}_\text{r}^*$--FM \qquad Black: $\widetilde{\text{u}}_\text{r}^*$--DFM}
\end{figure}

\subsection{A Physical Comparison}

A first inspection of radial velocity fitting and AT-INSAR inages, figures 4 - 7,
the gain with respect to the interferometric velocity $u_{ATI}$ may seem marginal.

To compare in terms of a physical quantity, we compute the relative error of associated kinetic energies. Results are shown in Table \ref{tab:RE_KE}.

\begin{table}[!h]
\begin{center}
{
\begin{tabular} {|c|c|c|c|}
\hline
\multicolumn{1}{|c|}{\textbf{Estimated solution}}  & \multicolumn{1}{|c|}{$\boldsymbol{\text{u}_\text{r}^*}$-NL} & \multicolumn{1}{c|} {$u_r^*$-FM} & \multicolumn{1}{c|} {$u_{ATI}$} \\
\hline
\textbf{RE of KE of the estimated solution} &  0.0582341 & 0.0219289 & 0.1180122 \\
\hline
\end{tabular}}
\end{center}
\caption[]{Relative errors of kinetic energies of the estimated solutions.}
\label{tab:RE_KE}
\end{table}

The relative error of the interferometric velocity is much greater that the fitted radial velocities. The error of using $u_{ATI}$ is about $11\%$.
In some applications this might be critical.

\subsection{Technical information}

For the functional versions of the Newton's methods and the surface simulations, we developed our own in house implementations.
The comparative performance of the methods above was carried out in a computer with the following specifications:

\subsubsection*{Hardware} \label{sec:HardwareSpecifications}
\begin{itemize}
\item Processor: AMD\,$^\copyright$ A10-5800B with Radeon(tm) HD Graphics, 3.80 GHz
\item Physical memory: 8.00 Gb (7.20 Gb usable).
\item Round-off unit (machine epsilon): $\epsilon_\mathcal{M} = 2.220446049250313 \times 10^{-16}$.
\end{itemize}

\subsubsection*{Software} \label{sec:SoftwareSpecifications}
\begin{itemize}
\item System type: 64-bit operating system.
\item Operating system: Windows 7 Professional\,$^\copyright$ 2009 Microsoft Corporation with Service Pack 1.
\item Programming language: Anaconda3 5.2.0 with Python 3.6.5 for 64 bits, Qt 5.9.4, PyQT5 5.9.2.
\item IDE: The Scientific PYthon Development EnviRonment (Spyder)\,$^\copyright$, version 3.2.8
\end{itemize}

For classical numerical methods, e.g., SVD decomposition and BFGS, we used the Python's routines.

\section{Conclusions}

Assuming the AT-INSAR-VB model, we have posed the radial velocities imaging problem, as the solution to a nonlinear integral equation. We have developed functional (infinite dimensional) versions of Modified Newton's methods, to solve this integral equation. Namely, a nonlinear system method coupled with Tikhonov regularization, and the BFGS method with known gradient for functional minimization.

For each technique, we have formulated the solution on function spaces, where the application of the Newton's method requires the Fr\'echet derivative of the objective functions Cheney 2001. 

Discrete models and numerical algorithms have been implemented. The numerical results are satisfactory. The functional approach leads to faster solutions in comparison with the classical discretize-then-optimize strategy. The fitting of the estimated radial velocity improves upon that of the $u_{ATI}$.  More over, the comparison of predicted Kinetic Energies, shows that in some applications, a better fit other than interferometric velocities is required.

This research is manifold, ocean waves modelling, sea surface imaging, computational methods, etc. On the modelling side, we have consider
only the swell spectrum.  It is of interest to consider for instance the JONSWAP and Pierson-Moskowitz spectra. 

Research on methods for sea surface imaging is ongoing. A straightforward computational continuation of this work, is the use of High Performance Computing. In our application,  an integral equation is solved for each point in the cross track coordinate. By the ATI-SAR-VB model, each solution is independent. Consequently, a parallel implementation in a low level computer language shall lead to even faster solutions.


\begin{thebibliography}{99}
\bibitem{Bao_Bruning_Alpers} M. Bao, C. Bruening, and W. Alpers. Simulation
of Ocean Waves Imaging by an Along-Track Interferometric Synthetic Aperture
Radar. IEEE Transactions on Geoscience and Remote Sensing, 35(3):618,631, May 1997.

\bibitem{Cheney} E. W. Cheney. Analysis for Applied Mathematics. Springer-Verlag, New
York, 2001.

\bibitem{Goldstein_Zebker} Goldstein, R. M., \& Zebker, H. A. (1987). Interferometric radar measurement of ocean surface currents. Nature, 328(6132), 707-709.

\bibitem{Hwang_etal} Hwang, P. A., Toporkov, J. V., Sletten, M. A., \& Menk, S. P. (2013). Mapping surface currents and waves with interferometric synthetic aperture radar in coastal waters: Observations of wave breaking in swell-dominant conditions. Journal of physical oceanography, 43(3), 563-582.

\bibitem{Mobley} C. D. Mobley. Modeling Sea Surfaces. A Tutorial on Fourier Transform
Techniques. Version 2.0. Sequoia Scientic, Inc., 2016.

\bibitem{Moreira_etal}A. Moreira, P. Prats-Iraola, M. Younis, G. Krieger, I. Hajnsek, and K. P. Papathanassiou.
A Tutorial on Synthetic Aperture Radar. IEEE Geoscience and Remote Sensing Magazine, 1(1):6, 43, March 2013.

\bibitem{Stuart} Stuart, A. M. (2010). Inverse problems: a Bayesian perspective. Acta numerica, 19, 451-559.

\bibitem{Zuazua} Zuazua, E. (2005). Propagation, observation, control and numerical approximation of waves. SIAM Review, 47(2), 197-243.
\end{thebibliography}
\end{document}